\documentclass{amsart}
\usepackage{tabu}
\usepackage{amssymb} 
\usepackage{amsmath} 
\usepackage{amscd}
\usepackage{amsbsy}
\usepackage{comment}
\usepackage[matrix,arrow]{xy}
\usepackage{hyperref}
\usepackage{color}
\usepackage{tabu}
\usepackage{longtable}
\usepackage{multirow}

\DeclareMathOperator{\Rad}{Rad}

\DeclareMathOperator{\Norm}{Norm}

\DeclareMathOperator{\Gal}{Gal}

\DeclareMathOperator{\ord}{ord}

\newcommand{\Q}{{\mathbb Q}}
\newcommand{\Z}{{\mathbb Z}}
\newcommand{\F}{{\mathbb F}}

\newcommand{\cE}{\mathcal{E}}

\newcommand{\cP}{\mathcal{P}}

\newcommand{\OO}{{\mathcal O}}

\newcommand{\sS}{\mathfrak{S}}

\newcommand{\fq}{\mathfrak{q}}

\begin {document}

\newtheorem{thm}{Theorem}
\newtheorem{lem}{Lemma}[section]

\theoremstyle{definition}

\theoremstyle{remark}

\title[]{On perfect powers that are sums of cubes of a three term arithmetic progression}

\author{Alejandro Arg\'{a}ez-Garc\'{i}a}
\address{Facultad de Ingenier\'{i}a Qu\'{i}mica, Universidad Aut\'{o}noma de Yucat\'{a}n. Perif\'{e}rico Norte Kil\'{o}metro 33.5, Tablaje Catastral 13615 Chuburna de Hidalgo Inn, M\'{e}rida, Yucat\'{a}n, M\'{e}xico. C.P. 97200 ; Centro de Invetigaci\'{o}n en Matem\'{a}ticas, A.C. - Unidad M\'{e}rida. Parque Cient\'{i}fico y Tecnol\'{o}gico de Yucat\'{a}n Km 5.5  Carretera Sierra Papacal - Chuburn\'{a} Puerto Sierra Papacal, M\'{e}rida, Yucat\'{a}n, M\'{e}xico. C.P. 97302}
\email{alejandro.argaez@correo.uady.mx}
\email{alejandro.argaez@cimat.mx}

\author{Vandita Patel}
\address{Department of Mathematics, University of Toronto, Bahen Centre, 40 St. George St., Room 6290, Toronto, Ontario, Canada, M5S 2E4}
\email{vandita@math.utoronto.ca}

\thanks{
The first-named author is supported by the Research Grant 221183 FOMIX, CONACYT - Government of Yucat\'an.
}

\date{\today}


\begin {abstract}
Using only elementary arguments, Cassels and Uchiyama (independently) 
determined all squares that are sums of three consecutive cubes.
Zhongfeng Zhang extended this result and determined all perfect powers that
are sums of three consecutive cubes. 
Recently, the equation $(x-r)^k + x^k + (x+r)^k$ has been studied for $k=4$ by Zhongfeng Zhang and for 
$k=2$ by Koutsianas.
In this paper, we complement
the work of Cassels, Koutsianas and Zhang by  considering the case when $k=3$ and showing that the equation 
$(x-r)^3+x^3+(x+r)^3=y^n$ with $n\geq 5$ a prime and $0 < r \leq 10^6$ only has trivial solutions $(x,y,n)$ which satisfy $xy=0$.
\end {abstract}
\maketitle

\section{Introduction} \label{intro}

In this paper, we prove the following: 
\begin{thm}\label{thm:main}
Let $p \ge 5$ be a prime. The only integer solutions to the equation
\begin{equation}\label{eq:main}
(x-r)^3 + x^3 + (x+r)^3 = y^p 
\qquad \gcd(x,r)=1
\end{equation}
with $0 < r\le 10^6$ are the trivial ones satisfying $xy=0$.
\end{thm}
The restriction $\gcd(x,r)=1$ is natural one, for 
otherwise it is easy to construct artificial solutions by
scaling. 
Our approach combines Frey curves and modularity with
a result of Mignotte based on linear form in logarithms
and various elementary techniques. This approach can easily
be extended to treat larger ranges of the parameter $r$.

Both our theorem and our methods extend theorems and ideas already found
in the literature. For example, 
Cassels \cite{Cassels} and Uchiyama \cite{Uchiyama} independently,
and using only elementary arguments, solved \eqref{eq:main}
with $p=2$ and $r=1$. 
Zhongfeng Zhang \cite{Zhang} extended the work of Cassels and Uchiyama and determined all perfect powers that
sums of three consecutive cubes; that is he solved \eqref{eq:main} with $r=1$ and  $p\geq 2$.

Power values of sums of powers of consecutive integers, or consecutive terms in
an arithmetic progression, have been considered by many authors from different
points of view. Techniques from algebraic number 
theory and Diophantine approximation have allowed the resolution of such equations with small exponents, as well as proofs of general
theorems. 
This can be seen in the work of  B\'{e}rces,  Pink,  Sava\c{s} and  Soydan \cite{BPSSoydan}, Brindza \cite{Br}, Cassels \cite{Cassels}, Gy\H{o}ry, Tijdeman and  Voorhoeve \cite{GTV},
Hajdu \cite{Ha}, Patel \cite{Patelsquares}, Patel and Siksek \cite{PatelSiksek}, Pint\'er \cite{P}, Sch\"affer \cite{Sc}, Schinzel and Tijdeman \cite{ScTi}, Soydan \cite{Soydan}, Urbanowicz \cite{U1}, 
and Zhang and Bai \cite{ZB} among many others.

Contributions from Bennett, Gy\H{o}ry and Pint\'er \cite{BGP}, Bennett Patel and Siksek \cite{BPS1}, \cite{BPS2} Pint\'er  \cite{P2} and Zhang \cite{Zhang},
resolved many such equations using Frey--Hellegouarch curves and Galois representations. 

Recently, the equations 
\[
(x-r)^k+x^k+(x+r)^k=y^n
\]
have been studied for $k=2$ by Koutsianas \cite{Koutsianas}, Koutsianas and Patel \cite{KoutsianasPatel} and for $k=4$ by Zhongfeng Zhang \cite{Zhang2}, via the study of Frey-Hellegouarch curves and Galois representations. 

\section{Small exponents: $p=2$ and $3$} \label{smallp}

For $p=2$, we have the following subfamily of equation~\eqref{eq:main}:
\[
(x-r)^3 + x^3 + (x+r)^3 = y^2.
\]
For a fixed value of $r$, 
since this defines an elliptic curve,
we may find all integer points via an appropriate computer package (for example \texttt{Magma}).

As $r$ varies, we immediately see that there are infinite families of solutions. One such family can be described by the following equations:
\begin{align*}
x &= 24a^2b^2, \\
r &= 9a^4 - 8b^4, \\
y &= 12ab(9a^4 - 8b^4),
\end{align*}
where $a$ and $b$ are integers. We easily see that $(x,r,y)$ are integral and satisfy the equation $(x-r)^3 + x^3 + (x+r)^3 = y^2$.

\medskip

For $p=3$, we consider the following subfamily of equation~\eqref{eq:main}:
\[
(x-r)^3 + x^3 + (x+r)^3 = y^3.
\]
We can define a smooth curve $C$ in $\mathbb{P}^2$:
\[
C \; : \; 3x(x^2+2r^2) = y^3.
\]
As this is a plane cubic, it is a curve of genus $1$. Let 
$P=(x,y,r)=(0,0,1) \in C(\mathbb{Q})$. Then $C$ is isomorphic
to an elliptic curve $E$ where the isomorphism $\phi \; : C \rightarrow E$
takes $P$ to the point at infinity on $E$. Using \texttt{Magma},
we find that $E$ has the model
\[
E \; : \; Y^2=X^3-648 \, .
\]
The map $\phi$ is given by $\phi(x:y:r)=(6y/x^3,-36r/x^2)$.
We find that $E(\Q) \cong \Z^2$ with Mordell--Weil basis
given by $(18,72)$ and $(9,9)$.

In particular, as $r$ varies, both cases $p=2$ and $p=3$ of equation~\eqref{eq:main} give rise to  infinitely many integer solutions.

\section{Four cases: equations of signature $(p,p,p)$ and $(p,p,2)$} \label{cases}

We may rewrite equation~\eqref{eq:main} as:
\[
3x(x^2+2r^2) = y^p \quad x,r,y,p \in \Z, \quad p\geq 2, \quad \gcd(x,r) = 1.
\]
Since $3 \mid y$ we can write $y=3w$. Thus
\begin{equation} \label{eq:mainfact}
x(x^2+2r^2)=3^{p-1} w^p.
\end{equation}
Note that $\gcd(x,x^2+2r^2)=1$ or $2$ according to
whether $2 \mid x$ or $2 \nmid x$. 
Thus, we consider four cases and apply a simple descent argument in each case
to obtain ternary equations.




\begin{itemize}
\item {\bf Case 1:} $2 \mid x$ and $3 \mid x$.
We rewrite equation~\eqref{eq:mainfact} as:
\[
x=2^{p-1} \cdot 3^{p-1} w_1^p, \qquad x^2+2r^2=2 w_2^p.
\]
Substituting and simplifying we obtain:
\begin{equation}\label{eqn:case1}
w_2^p-2^{2p-3} \cdot 3^{2p-2} w_1^{2p}=r^2.
\end{equation}
\item {\bf Case 2:} $2 \mid x$ and $3 \nmid x$.
We rewrite equation~\eqref{eq:mainfact} as:
\[
x=2^{p-1} w_1^p, \qquad x^2+2r^2= 2 \cdot 3^{p-1} w_2^p.
\]
Substituting and simplifying we obtain:
\begin{equation}\label{eqn:case2}
3^{p-1} w_2^p-2^{2p-3} w_1^{2p}=r^2.
\end{equation}
\item {\bf Case 3:} $2 \nmid x$ and $3 \mid x$.
We rewrite equation~\eqref{eq:mainfact} as:
\[
x=3^{p-1} w_1^p, \qquad x^2+2r^2= w_2^p.
\]
Substituting and simplifying we obtain:
\begin{equation}\label{eqn:case3}
w_2^p-3^{2p-2} w_1^{2p}=2r^2.
\end{equation}
\item {\bf Case 4:} $2 \nmid x$ and $3 \nmid x$.
We rewrite equation~\eqref{eq:mainfact} as:
\[
x=w_1^p, \qquad x^2+2r^2= 3^{p-1} w_2^p.
\]
Substituting and simplifying we obtain:
\begin{equation}\label{eqn:case4}
3^{p-1} w_2^p-w_1^{2p}=2r^2.
\end{equation}

\end{itemize}

Cases 1 and 2 give rise to ternary equations of signature $(p,p,2)$ in the variables $(x,r,y,p)$. 
Here, we are able to prove Theorem~\ref{thm:main} without any restriction on $r$ 
using Chabauty's method for $p=5$ and the modular approach for $p \geq 7$. 
This is expanded upon in Sections~\ref{sec:p=5} and \ref{sec:case12} respectively.

Cases 3 and 4 also produce equations of signature $(p,p,2)$ in the variables 
$(x,r,y,p)$, however, these equations currently don't seem tractable under the modular method. 
Therefore, we consider their ternary equation as having signature $(p,p,p)$ in the variables $(x,y,p)$ 
when fixing values of $r$. Moving away from the modular method and taking a more classical approach yields 
results for $r \leq 10^6$ with minimum computational effort. This is expanded upon in 
the remaining sections. 
\section{The case $p=5$}
\label{sec:p=5}

For $p=5$, 
equations \eqref{eqn:case1}, \eqref{eqn:case2}, \eqref{eqn:case3}, \eqref{eqn:case4} can 
be rewritten as curves of genus $2$ defined over $\Q$. 
For a curve $C$ of genus $2$ defined over $\Q$, the method of Chabauty and Coleman
can be used to determine the rational points provided the rank of the Mordell--Weil
group of the Jacobian is at most $1$. For an introduction to the method 
of Chabauty and Coleman we recommend \cite{MP} and \cite{SiksekChab}.
Fortunately the method is implemented in the \texttt{Magma} software
package and we were able to solve equations \eqref{eqn:case1}, \eqref{eqn:case2}
and \eqref{eqn:case3} for $p=5$ using this implementation.

\subsection{Case 1}
Let $p=5$ and consider \eqref{eqn:case1}. Write
$X=w_2/w_1^2$ and $Y=r/w_1^5$. Then $(X,Y)$ is a rational
point on the genus $2$ curve
\[
C \; : \; Y^2=X^5-2^7 \cdot 3^8.
\]
Let $J$ be the Jacobian on $C$. Using \texttt{Magma}, which uses
a $2$-descent algorithm of Stoll \cite{Stoll}, we find that the Mordell--Weil
rank of $J$ is $0$. The \texttt{Magma} implementation of Chabauty gives
 $C(\Q)=\{\infty\}$. Thus $w_1=0$ and $w_2=1$, which gives the 
trivial solution to \eqref{eq:main}.

\subsection{Case 2}
We now consider \eqref{eqn:case2} with $p=5$. With $X$, $Y$ as above we see
that $(X,Y)$ is a rational point on 
\[
C \; : \; Y^2=3^4 X^5-2^7.
\]
The Mordell--Weil group of the Jacobian is infinite cyclic, 
generated by the point, given in Mumford representation,
\[
\left(X^2 + \frac{20}{81} X + \frac{8}{9}, \; \frac{1124}{81} X + \frac{32}{9} \right).
\]
The \texttt{Magma} Chabauty implementation again gives $C(\Q)=\{\infty\}$
and this gives the trivial solution to \eqref{eq:main}.

\subsection{Case 3}
Let $p=5$ and consider \eqref{eqn:case3}.
We let $X=w_2/w_1^2$ and $Y=2r/w_1^5$. Then $(X,Y)$ is a rational point
on 
\[ 
C \; : \; Y^2=2(X^5-3^8).
\]
Let $P=(9,324) \in C(\Q)$. The Mordell--Weil group of the Jacobian
is infinite cyclic, generated by $P-\infty$. The Chabauty implementation
now gives
$\infty$, $(9,324)$, $(9,-324)$ as the only rational points. Thus
$w_1=0$, $w_2=1$ or $w_1=1$, $w_2=9$, $r=\pm 162$. The former gives a
trivial solution to \eqref{eq:main} and the latter can be discarded
as it does not
satisfy $\gcd(x,r)=1$.

\subsection{Case 4}
We mention in passing that for $p=5$ the Mordell--Weil rank of the Jacobian
of the genus $2$ curve corresponding to equation~\eqref{eqn:case4} is $2$, and so
Chabauty's method is inapplicable.

\section{Cases 1 and 2: a modular approach} \label{sec:case12}

In this section, we are able to give a complete resolution to Cases 1 and 2 for any
$r \in \Z$ and $p \ge 7$. 
In order to achieve this,  we make use of the work of Bennett and Skinner 
\cite{BenS} (see also \cite{Siksek}). That work
builds upon the celebrated modularity theorem of Wiles, Breuil, Conrad,
Diamond and Taylor \cite{Wiles}, \cite{modularity}, on Ribet's level lowering theorem
\cite{Ribet}, and on Mazur's isogeny theorem \cite{Mazur}.

\begin{lem}\label{lem:cases12mod}
Let $p\geq 7$ be a prime. Suppose $2 \mid x$. Then equation~\eqref{eq:main} has no integer solutions.
\end{lem}
The reader will note that the condition $2 \mid x$ corresponds to Cases 1 and 2 above.
\begin{proof}
Suppose first that $2 \mid x$ and $3 \mid x$ and thus we are in Case 1. Thus we have a
non-trivial solution to a ternary equation of signature $(p,p,2)$, namely \eqref{eqn:case1}.
Our condition $\gcd(x,r)=1$ readily implies that the three terms in \eqref{eqn:case1}
are coprime.
Changing the sign of $r$ if necessary, we may suppose that
$r \equiv 1 \pmod{4}$.
Using the recipes of Bennett and Skinner \cite{BenS}
we attach to this the Frey-Hellegouarch curve:
\[
E \; : \; Y^2+XY=X^3+\left(\frac{r-1}{4} \right) X^2-2^{2p-9} \cdot 3^{2p-4} w_1^{2p} X.
\]
By \cite[Lemma 3.3]{BenS}, 
$E$ arises modulo $p$ from a newform of weight $2$ and level $6$.
As there are no such newforms we have a contradiction.


Hence, we suppose that $2 \mid x$ and $3 \nmid x$ (Case 2). We now consider the 
ternary equation \eqref{eqn:case2} where again the three terms are coprime
and again we may suppose $r \equiv 1 \pmod{4}$.
%
The recipes
of Bennett and Skinner \cite{BenS} attach to this the following Frey-Hellegouarch curve:
\[
E \; : \; Y^2+XY=X^3+\left( \frac{r-1}{4}\right) X^2-2^{2p-9} w_1^{2p} X.
\]
This too arises from a newform of weight $2$ and level $6$, and again, we arrive at a contradiction.
\end{proof}

\section{Cases 3 and 4: Bounding $p$} \label{Boundp}


The aim of this section is to bound $p$ in Cases 3 and 4. 
We require the following theorem of Mignotte, which is found in \cite[Chapter~12, p.~423]{Cohen2}.

\begin{thm}[Mignotte]\label{thm:Mignotte}
Assume that the exponential Diophantine inequality 
\[
|ax^n - by^n | \leq c, \quad \text{ with } a,b,c \in \Z_{\ge 0} \text{ and }  a\neq b,
\]
has a solution in strictly positive integers $x$ and $y$ with $\max\{x,y\} > 1$. Let $A = \max \{a,b,3\}$. Then
\[
n \leq \max \left\{ 3 \log\left(1.5| c/b| \right),\; \frac{7400\log A}{\log\left(1+\log A/\log\left(\mid a/b\mid\right)\right)}\right\}.
\]
\end{thm}

A simple application of Mignotte's Theorem yields workable bounds for
$p$ for equations~\eqref{eqn:case3} and
\eqref{eqn:case4}, and hence for Cases 3 and 4 respectively. 
We first deal with Case 3.  
\begin{lem}
Suppose $2 \nmid x$ and $3 \mid x$ and $ |r| \leq 1.55\times 10^{1697}$. 
If equation~\eqref{eq:main} has an integral non-trivial solution, then $p<24,000$.
\end{lem}
\begin{proof}
Multiplying \eqref{eqn:case3} by $9$ gives
\begin{equation}\label{eqn:Linlog3}
9w_2^p - (3 w_1)^{2p}=18r^2.
\end{equation}
We apply Theorem~\ref{thm:Mignotte} to equation~\eqref{eqn:Linlog3} with $a=9$, $b=1$, $c=18r^2$ and we obtain
\begin{align*}
p &\leq \max \left\{ 3 \log\left(3^3 \cdot r^2  \right),\; 7400\log 9 / \log 2 \right\}\\
 & \leq \max \left\{ 23457.4384\ldots ,\; 23457.4450\ldots \right\}.
\end{align*}
We deduce that $p < 24,000$, which completes the proof.
\end{proof}

Next we deal with Case 4.
\begin{lem}
Suppose $2 \nmid x$ and $3 \nmid x$. Suppose $ |r| \leq 2.99\times 10^{848}$. 
If equation~\eqref{eq:main} has an integral non-trivial solution, then $p\leq 12,000$.
\end{lem}
\begin{proof}
Multiplying \eqref{eqn:case4} by $3$ gives
\begin{equation}\label{eqn:Linlog4}
(3w_2)^p- 3\cdot w_1^{2p}=6r^2.
\end{equation}
We apply Theorem~\ref{thm:Mignotte} to equation~\eqref{eqn:Linlog4} with $a=3$, $b=1$, $c=6r^2$ 
to obtain
\begin{align*}
p &\leq \max \left\{ 3 \log\left( 3^2\cdot r^2  \right),\; 7400\log 3 / \log 2 \right\}\\
 & \leq \max \left\{ 11727.7611\ldots ,\; 11728.7225\ldots \right\}.
\end{align*}
Hence,if equation~\eqref{eq:main} has an integral non-trivial solution, then $p < 12,000$.
\end{proof}

\section{Criteria for eliminating equations of signature $(p,2p,2)$} \label{sec:criteria}

To prove Theorem~\ref{thm:main} we would like to solve equation \eqref{eq:main}
for $0 < r \le 10^6$. We have seen that it is enough to solve equations \eqref{eqn:case3}
and \eqref{eqn:case4} for $0<r \le 10^6$. In the previous section we saw that
the assumption $0<r\leq 10^6$ forces 
$p <24000$ in \eqref{eqn:case3} and $p<12000$ in \eqref{eqn:case4}.
This leaves us with a finite, albeit very large, number of equations
to solve. These all have form
\begin{equation}\label{eqn:tern}
ay_2^{p} - by_1^{2p} = c
\end{equation}
where $p$ is a fixed odd prime and $a$, $b$ and $c$ are positive integers
satisfying $\gcd(a,b,c)=1$.
In this section, we state various criteria for equations of the form \eqref{eqn:tern}
to have no integer solutions.
Proofs of these criteria can be found  in \cite{BPS2} and \cite{Patelthesis}; the latter providing extensive details
and motivation.  
In later sections, we subject equations \eqref{eqn:case3} and \eqref{eqn:case4} 
to these criteria.
 
\subsection{Criteria for the non-existence of solutions} \label{SG}

\begin{lem}\label{lem:Sophiecriterion}
Let $p \ge 3$ be a prime. Let $a$, $b$ and $c$ be positive integers such that 
$\gcd(a,b,c)=1$. Let $q=2k p+1$ be a prime that does
not divide $a$. Define
\begin{equation}\label{eqn:mu}
\mu(p,q)=\{ \eta^{2p} \; : \; \eta \in \F_q \}
=\{0\} \cup \{ \zeta \in \F_q^* \; : \; \zeta^{k}=1\}
\end{equation}
and
$$
B(p,q)=\left\{ \zeta \in \mu(p,q) \; : \; ((b \zeta+c)/a)^{2k} \in \{0,1\} \right\} \, .
$$
If $B(p,q)=\emptyset$,
then equation~\eqref{eqn:tern} 
does not have integral solutions.
\end{lem}

The proof is straightforward and can be found in \cite{BPS2} and \cite{Patelthesis}.

\subsection{Local Solubility}\label{sub:locsol}
In this subsection, we appeal to classical local solubility methods 
for Diophantine equations in order  
to conclude non-existence of solutions for many tuples $(a,b,c,p)$
in equation~\eqref{eqn:tern}.
 
In \eqref{eqn:tern}, recall the condition  $\gcd(a,b,c)=1$.
Let $g=\Rad(\gcd(a,c))$ and suppose that $g > 1$. Then $g \mid y_1$,
and we can write $y_1=g y_1^\prime$. Thus
\[
a y_2^p- b g^{2 p} {y_1^\prime}^{2p}=c. 
\]
Removing a factor of $g$ from the coefficients, 
we obtain
\[
a^\prime {y_2}^p - b^\prime {y_1^\prime}^{2p}=c^\prime,
\]
where $a^\prime=a/c$ and $c^\prime=c/g<c$. Similarly, if $h=\gcd(b,c)>1$, we obtain
\[
a^\prime {y_2^\prime}^p - b^\prime {y_1}^{2p}=c^\prime,
\]
where $c^\prime=c/h<c$. Applying these operations repeatedly, 
we arrive at an equation of the form
\begin{equation}\label{eqn:predescent}
A \rho^p-B \sigma^{2p}=C
\end{equation}
where $A$, $B$, $C$ are now pairwise coprime.
A necessary condition for the existence of solutions is that
for any odd prime $q \mid A$, the residue $-BC$ modulo $q$ is a square.
Besides this basic test, we also check for local solubility
at the primes dividing $A$, $B$, $C$, and all primes $q \le 19$.


\subsection{A Descent}\label{sub:furtherdesc}

If local techniques of Section~\ref{sub:locsol} fail to rule out solutions 
to equation \eqref{eqn:tern} for particular coefficients and exponent $(a,b,c,p)$ 
then we may perform a further descent to rule out solutions. 
With $A$, $B$, $C$ as in \eqref{eqn:predescent}
we let
\[
B^\prime=\prod_{\text{$\ord_q(B)$ is odd}} q.
\]
Thus $B B^\prime=v^2$. Write $A B^\prime=u$ and $C B^\prime=m n^2$
with $m$ squarefree. Rewrite \eqref{eqn:predescent}
as
\[
(v \sigma^p+n \sqrt{-m})(v \sigma^p-n \sqrt{-m})=u \rho^p. 
\]
Let $K=\Q(\sqrt{-m})$ and $\OO$ be its ring of integers.
Let $\sS$ contain the prime ideals of $\OO$ that
divide $u$ or $2n \sqrt{-m}$. 
Clearly 
$(v\sigma^p+n \sqrt{-m}) {K^*}^p$ belongs to the ``$p$-Selmer group''
\[
K(\sS,p)=\{\epsilon \in K^*/{K^*}^p \; : \; 
\text{$\ord_\cP(\epsilon) \equiv 0 \mod{p}$ for all $\cP \notin \sS$}
\}.
\]
This is an $\F_p$-vector space of finite dimension 
(\cite[Proof of Proposition VIII.1.6]{Silverman}), and
can be computed by \texttt{Magma} using the command \texttt{pSelmerGroup}.
Let
\[
\cE=\{ \epsilon \in K(\sS,p) \; : \; \Norm(\epsilon)/u \in {\Q^*}^p \}.
\]
It follows that
\begin{equation}\label{eqn:furtherdescent}
v \sigma^p+n \sqrt{-m}=\epsilon \eta^p,
\end{equation}
where $\eta \in K^*$ and $\epsilon \in \cE$.

\begin{lem}\label{lem:valuative}
Let $\fq$ be a prime ideal of $K$. Suppose one of the following
holds:
\begin{enumerate}
\item[(i)] $\ord_\fq(v)$, $\ord_\fq(n\sqrt{-m})$, $\ord_\fq(\epsilon)$
are pairwise distinct modulo $p$;
\item[(ii)] $\ord_\fq(2v)$, 
$\ord_\fq(\epsilon)$, $\ord_\fq(\overline{\epsilon})$
are pairwise distinct modulo $p$;
\item[(iii)] $\ord_\fq(2 n \sqrt{-m})$, 
$\ord_\fq(\epsilon)$, $\ord_\fq(\overline{\epsilon})$
are pairwise distinct modulo $p$.
\end{enumerate}
Then there is no $\sigma \in \Z$
and $\eta \in K$ satisfying \eqref{eqn:furtherdescent}.
\end{lem}

Proof and further details may be found in \cite{BPS2} and \cite{Patelthesis}.

\begin{lem}\label{lem:furtherdescent}
Let $q=2k p+1$ be a prime. 
Suppose $q\OO=\fq_1 \fq_2$ where $\fq_1$, $\fq_2$
are distinct, and such that $\ord_{\fq_j}(\epsilon)=0$
for $j=1$, $2$. Let 
\[
\chi(p,q)=\{ \eta^p \; : \; \eta \in \F_q \}.
\]
Let
\[
C(p,q)=\{\zeta \in \chi(p,q) \; : \;
((v \zeta+n\sqrt{-m})/\epsilon)^{2k} \equiv \text{$0$ or $1 \mod{\fq_j}$
for $j=1$, $2$}\}.
\]
Suppose $C(p,q)=\emptyset$. 
Then there is no $\sigma \in \Z$
and $\eta \in K$ satisfying \eqref{eqn:furtherdescent}.
\end{lem}

The proof of Lemma~\ref{lem:furtherdescent} is a simple modification 
of the proof of Lemma~\ref{lem:Sophiecriterion}, the latter can be found in \cite{BPS2} and \cite{Patelthesis}.

\section{Application to Case 3}\label{sec:Appcase3}

We wrote a \texttt{Magma} program that subjected equation~\eqref{eqn:case3} with $7 \le p<24,000$ 
and $0 <r\leq 10^6$
to Lemma~\ref{lem:Sophiecriterion}, local solubility from subsection~\ref{sub:locsol} as well as lemmata~\ref{lem:valuative} and \ref{lem:furtherdescent}. 
Recall the assumption that $2 \nmid x,\; 3 \mid x$ and $\gcd(x,r)=1$ for Case 3. Thus $3 \nmid r$.  This allowed us to reduce the number of values of $r$ that we subject to the tests of Section~\ref{sec:criteria}.
The full time of computations took roughly $66$ hours.
The results are summarised in Table~\ref{table:3}.

\begin{table}[!htbp]
\centering
{
\begin{tabu}{|c|c|c|c|}
\hline
\multirow{3}{*}{Exponent $p$} & equations~\eqref{eqn:case3} surviving  & number surviving  & number surviving\\
				 &  Lemma~\ref{lem:Sophiecriterion} &   after local & after further\\
				 &  with exponent $p$ & solubility tests  & descent\\
\hline\hline
$7$ & $67,600$  &  $29,201$ & $4$ \\
\hline
$11$ & $7,233$ &  $5,334$ &  $2$\\
\hline
$13$ & $1,718$  &  $808$ & $2$  \\
\hline
$17$ & $848$  &   $453$ & $1$  \\
\hline
$19$ & $2,032$  & $1,023$ &  $1$  \\
\hline
$23$ & $19$ &  $12$ &  $1$\\
\hline
$29$ & $8$  &  $5$ & $1$  \\
\hline
$31$ & $47$  &   $24$ & $1$  \\
\hline
$37$ & $5$  & $5$ &  $1$  \\
\hline
$41 \le p <24,000$ & $0$ & $0$ & $0$\\
\hline\hline
Total & $79,510$ & $36,865$ & $14$ \\
\hline\hline
\end{tabu}
}
\vskip1ex
\caption{}
\label{table:3}
\end{table}

\section{Application to Case 4}\label{sec:Appcase4}

We subjected equation~\eqref{eqn:case4} with $5 \le p<12,000$ and $0 < r\leq 10^6$ to 
Lemma~\ref{lem:Sophiecriterion}, local solubility from Section~\ref{sub:locsol} 
as well as lemmata~\ref{lem:valuative} and \ref{lem:furtherdescent}. 
Equation~\eqref{eqn:case4} is under the assumption that $2,3 \nmid x$ and $\gcd(x,r)=1$.
Since $x^2 + 2r^2 = 3^{p-1}w_1^p$, we easily conclude that $3 \nmid r$.  
This allowed us to significantly cut down the number of values of $r$ that we subject to the tests 
of Section~\ref{sec:criteria}. 
The full time of computations took roughly $70$ hours. 
The results are summarised in Table~\ref{table:4}.

\begin{table}[!htbp]
\centering
{
\begin{tabu}{|c|c|c|c|}
\hline
\multirow{3}{*}{Exponent $p$} & equations~\eqref{eqn:case4} surviving  & number surviving  & number surviving\\
				 &  Lemma~\ref{lem:Sophiecriterion} &   after local & after further\\
				 &  with exponent $p$ & solubility tests  & descent\\
\hline\hline
$5$ & $219,921$ &  $65,821$ & $777$ \\
\hline
$7$ & $25,308$  &  $12,994$ & $0$ \\
\hline
$11$ & $3,201$ &  $1,834$ &  $0$\\
\hline
$13$ & $948$  &  $439$ & $0$  \\
\hline
$17$ & $385$  &   $201$ & $0$  \\
\hline
$19$ & $1,825$  & $1,108$ &  $0$  \\
\hline
$23$ & $7$ &  $3$ &  $0$\\
\hline
$29$ & $1$  &  $1$ & $0$  \\
\hline
$31$ & $3$  &   $2$ & $0$  \\
\hline
$37$ & $3$  & $2$ &  $0$  \\
\hline
$41$ & $1$  & $1$ &  $0$  \\
\hline
$43 \le p <12,000$ & $0$ & $0$ & $0$\\
\hline\hline
Total & $251,603$ & $82,406$ & $777$ \\
\hline\hline
\end{tabu}
}
\vskip1ex
\caption{}
\label{table:4}
\end{table}

\section{Modularity, Thue equations and completing the proof} \label{Thue}

We notice that in Table~\ref{table:3}, equations remain to be solved with an unusually high exponent $p$, and thus Magma's Thue Solver will not yield results. 

On closer inspection, one notices that the tuples $(a,b,c)$ that remain (which describe equation~\eqref{eqn:tern}) are of the form $(a,b,2^k)$ with $k \geq 5$. Hence, these cases may be eliminated by once again, appealing to the modular method.

\begin{lem}\label{lem:modpow2}
Suppose $p\geq 5$ is a prime. Suppose $2 \nmid x$ and $3 \mid x$. Suppose $r = 2^k$ with $k \geq 2$. Then equation~\eqref{eq:main} has no integer solutions.
\end{lem}

\begin{proof}
Since $2 \nmid x$ and $3 \mid x$, recall that we obtain the ternary equation~\eqref{eqn:case3} of signature $(p,p,p)$:
\[
w_2^p-3^{2p-2} w_1^{2p}=2^{2k+1},
\]
with all three terms coprime (this follows readily from our imposed condition $\gcd(x,r)=1$).
Following Kraus \cite{KrausFermat}, we attach to this the Frey-Hellegouarch curve
\[
E \; : \; Y^2 = X(X - 3^{2p-2}w_1^{2p} )(X + 2^{2k+1}).
\]
By \cite[Section 4]{KrausFermat}  
(which again builds on the aforementioned deep theorems of Mazur, Ribet,
Wiles, Conrad, Diamond and Taylor)
this arises modulo $p$ from a newform of weight $2$ and level $6$. 
As there are no such newforms we arrive at a contradiction.
\end{proof}

Lemma~\ref{lem:modpow2} allows us to eliminate $13$ equations from Table~\ref{table:3}, all with large exponent $p$. Thus, 
for Case 3 it remains to solve only $1$ equation with exponent $p=7$.

\subsection{Solving the final equations}
To complete the proof of Theorem~\ref{thm:main} it remains to solve $1+777=778$ equations of the form \eqref{eqn:tern}. Letting 
$\sigma = w_2$ and $\tau = w_1^2$,
these can be considered as Thue equations,
\[
a\sigma^p - b \tau^p = c \, .
\]
where the exponent $p=5$ in all but one case where $p=7$.
Finally, we appeal to \texttt{Magma's} Thue solver \cite{magma}; details of the algorithms employed by this 
solver can be found in \cite[Chapter VII]{Smart}. Solving these Thue equations we have found no non-trivial
solutions.
This completes the proof of Theorem~\ref{thm:main}.



\end{document}